# Properties of Soft Semi-open and Soft semi-closed Sets

Sabir Hussain

Department of Mathematics, College of Science, Qassim University

P. O. box 6644, Buraydah 51452, Saudi Arabia

E-mail: sabiriub@yahoo.com; sh.hussain@qu.edu.sa

**Abstract**
Molodstov[10] introduced soft set theory as a new mathematical approach for solving problems having uncertainties. Many researchers worked on the findings of structures of soft set theory and applied to many problems having uncertainties. Recently, Bin Chen [3-4] introduced and explored the properties of soft semi-open sets and soft-semi-closed sets in soft topological spaces. In this paper we continue to investigate the properties of soft semi-open and soft semi-closed sets in soft topological spaces. We define soft semi-exterior, soft semi-boundary, soft semi-open neighbourhood and soft semi-open neighbourhood systems in soft topological spaces. Moreover we discuss the characterizations and properties of soft semi-interior, soft semi-exterior, soft semi-closure and soft semi-boundary. We also develop the relationship between soft semi-clopen sets and soft semi-boundary. The addition of this topic in literature will strengthen the theoretical base for further applications of soft topology in decision analysis and information systems.

**1. Introduction**
In 1999, Molodtsov [10] introduced soft sets theory as a new general mathematical approach to deal with uncertainties and unclear defined objects. In [11], they applied successfully this approach for modelling the problems in engineering physics, computer science, economics, social sciences and medical sciences in directions such as, smoothness of functions, game theory, operations research, Riemann-integration, Perron integration, probability and theory of measurement.

      In soft systems a very general frame work has been provided with the involvement of parameters. Therefore the researches work on soft set theory and its applications in various disciplines and real life problem are now catching momentum. Maji et. al [8] defined and studied several basic notions of soft set theory and discussed in detail the application of soft set theory in decision making problems [9]. Xiao et. al [15] and Pei et. al [13] studied the relationship between soft sets and information systems. Kostek [6] presented the criteria to measure sound quality using approach of soft sets. Mushrif et. al [12] used the notions of soft set theory to develops the remarkable method for the classification of natural textures. Many researchers contributed towards the algebraic structures of soft set theory.

      In 2011, Shabir and Naz [14] initiated the study of soft topological spaces and defined basic notions of soft topological spaces. After that, Hussain et. al [5] [1], Aygunoglu et.al [2], Zoultrana et. al [16] added many notions and concepts towards the properties of soft topological spaces. Recently in 2013, Bin Chen [3-4] introduced and explored the properties of soft semi-open sets and soft-semi-closed sets in soft topological spaces.

      In this paper we continue to investigate the properties of soft semi-open sets and soft semi-closed sets in soft topological spaces. We define soft semi-exterior, soft semi-boundary, soft semi-open neighborhood and soft semi-open neighborhood systems in soft topological spaces. Moreover we discuss the characterizations and properties of soft semi-interior, soft semi-exterior, soft semi-closure and soft semi-boundary in soft topological spaces.

---

## 2. Preliminaries

First we recall some definitions and results defined and discussed in [1, 3, 4, 10, 11, 14, 16], which will use in the sequel.

**Definition 2.1.** Let $X$ be an initial universe and $E$ be a set of parameters. Let $P(X)$ denotes the power set of $X$ and $A$ be a non-empty subset of $E$. A pair $(F, A)$ is called a soft set over $X$, where $F$ is a mapping given by $F: A \to P(X)$. In other words, a soft set over $X$ is a parameterized family of subsets of the universe $X$. For $e \in A$, $F(e)$ may be considered as the set of $e$-approximate elements of the soft set $(F, A)$. Clearly, a soft set is not a set.

**Definition 2.2.** For two soft sets $(F, A)$ and $(G, B)$ over a common universe $X$, we say that $(F, A)$ is a soft subset of $(G, B)$ if

(1) $A \subseteq B$ and

(2) for all $e \in A$, $F(e)$ and $G(e)$ are identical approximations. We write $(F, A) \widetilde{\subseteq} (G, B)$.

$(F, A)$ is said to be a soft super set of $(G, B)$, if $(G, B)$ is a soft subset of $(F, A)$. We denote it by $(F, A) \widetilde{\supseteq} (G, B)$.

**Definition 2.3.** Two soft sets $(F, A)$ and $(G, B)$ over a common universe $X$ are said to be soft equal, if $(F, A)$ is a soft subset of $(G, B)$ and $(G, B)$ is a soft subset of $(F, A)$.

**Definition 2.4.** The union of two soft sets of $(F, A)$ and $(G, B)$ over the common universe $X$ is the soft set $(H, C)$, where $C = A \cup B$ and for all $e \in C$,

$$H(e) = \begin{cases} F(e), & \text{if } e \in A - B \\ G(e), & \text{if } e \in B - A \\ F(e) \cup G(e), & \text{if } e \in A \cap B \end{cases}$$

We write $(F, A) \widetilde{\cup} (G, B) = (H, C)$.

**Definition 2.5.** The intersection $(H, C)$ of two soft sets $(F, A)$ and $(G, B)$ over a common universe $X$, denoted $(F, A) \widetilde{\cap} (G, B)$, is defined as $C = A \cap B$, and $H(e) = F(e) \widetilde{\cap} G(e)$, for all $e \in C$.

**Definition 2.6.** The difference $(H, E)$ of two soft sets $(F, E)$ and $(G, E)$ over $X$, denoted by $(F, E) \widetilde{\setminus} (G, E)$, is defined as $H(e) = F(e) \setminus G(e)$, for all $e \in E$.

**Definition 2.7.** Let $(F, E)$ be a soft set over $X$ and $Y$ be a non-empty subset of $X$. Then the sub soft set of $(F, E)$ over $Y$ denoted by $(Y_F, E)$, is defined as follows: $F_Y(\alpha) = Y \widetilde{\cap} F(\alpha)$, for all $\alpha \in E$. In other words $(Y_F, E) = \widetilde{Y} \widetilde{\cap} (F, E)$.

**Definition 2.8.** The relative complement of a soft set $(F, A)$ is denoted by $(F, A)'$ and is defined by $(F, A)' = (F', A)$ where $F': A \to P(U)$ is a mapping given by $F'(\alpha) = U \setminus F(\alpha)$, for all $\alpha \in A$.

**Definition 2.9.** Let $x \in X$, then $(x, E)$ denotes the soft set over $X$ for which $x(\alpha) = \{x\}$, for all $\alpha \in E$.

**Definition 2.10.** Let $(F, E)$ be a soft set over $X$ and $Y$ be a non-empty subset of $X$. Then the sub soft set of $(F, E)$ over $Y$ denoted by $(Y_F, E)$, is defined as follows: $F_Y(\alpha) = Y \cap F(\alpha)$, for all

$\alpha \in E$. In other words $(Y_F, E) = \tilde{Y} \cap (F, E)$.

**Definition 2.11.** Let $(F, E)$ be a soft set over $X$ and $x \in X$. We say that $x \in (F, E)$ read as $x$ belongs to the soft set $(F, E)$, whenever $x \in F(\alpha)$, for all $\alpha \in E$. Note that $x \in X$, $x \notin (F, E)$, if $x \notin F(\alpha)$ for some $\alpha \in E$.

**Definition 2.12.** Let $\tau$ be the collection of soft sets over $X$, then $\tau$ is said to be a soft topology on $X$ if

(1) $\Phi$, $\tilde{X}$ belong to $\tau$.
(2) the union of any number of soft sets in $\tau$ belongs to $\tau$.
(3) the intersection of any two soft sets in $\tau$ belongs to $\tau$.

The triplet $(X, \tau, E)$ is called a soft topological space over $X$.

**Definition 2.13.** Let $(X, \tau, E)$ be a soft topological space over $X$ then soft interior of soft set $(F, E)$ over $X$ is denoted by $(F, E)°$ and is defined as the union of all soft open sets contained in $(F, E)$. Thus $(F, E)°$ is the largest soft open set contained in $(F, E)$. A soft set $(F, E)$ over $X$ is said to be a soft closed set in $X$, if its relative complement $(F, E)'$ belongs to $\tau$.

**Definition 2.14.** Let $(X, \tau, E)$ be a soft topological space over $X$ and $(F, E)$ be a soft set over $X$. Then the soft closure of $(F, E)$, denoted by $\overline{(F, E)}$ is the intersection of all soft closed super sets of $(F, E)$. Clearly $\overline{(F, E)}$ is the smallest soft closed set over X which contains $(F, E)$.

**Definition 2.15.** Let $(X, \tau, E)$ be a soft topological space over $X$ and $(F, E)$ be a soft set over $X$. Then $(F, E)$ is called soft semi-open set if and only if there exists a soft open set $(G, E)$ such that $(G, E) \tilde{\subseteq} (F, E) \tilde{\subseteq} \overline{(F, E)}$. The set of all soft semi-open sets is denoted by $S.S.O(X)$. Note that every soft open set is soft semi-open set.

A sot set $(F, E)$ is said to be soft semi-closed if its relative complement is soft semi-open. Equivalently there exists a soft closed set $(G, E)$ such that $(G, E)° \tilde{\subseteq} (F, E) \tilde{\subseteq} (G, E)$. Note that every soft closed set is soft semi-closed set.

**Definition 2.16.** Let $(X, \tau, E)$ be a soft topological space over $X$.

(i) soft semi-interior of soft set $(F, E)$ over $X$ is denoted by $int^s(F, E)$ and is defined as the union of all soft semi-open sets contained in $(F, E)$.

(ii) soft closure of $(F, E)$ over $X$ is denoted by $cl^s(F, E)$ is the intersection of all soft semi-closed super sets of $(F, E)$.

### 3. Properties of soft semi-open sets and soft semi-closed sets

Hereafter, we will denote the soft semi-interior, soft semi-exterior, soft semi-closure and soft semi-boundary as $int^s$, $ext^s$, $cl^s$ and $bd^s$ respectively.

**Definition 3.1.** Let $(X, \tau, E)$ be a soft topological space over $X$, $(G, E)$ be a soft set over $X$ and $x \in X$. Then $x$ is said to be a soft semi-interior point of $(G, E)$, if there exists a soft semi-open set $(F, E)$ such that $x \in (F, E) \tilde{\subseteq} (G, E)$.

**Definition 3.2.** Let $(X, \tau, E)$ be a soft topological space over $X$, $(G, E)$ be a soft set over $X$ and $x \in X$. Then $(G, E)$ is said to be a soft semi-neighborhood of $x$, if there exists a soft semi-open set $(F, E)$ such that $x \in (F, E) \tilde{\subseteq} (G, E)$.

**Proposition 3.3.** Let $(X, \tau, E)$ be a soft topological space over $X$. For any soft semi-open set $(F, E)$

over $X$, $(F,E)$ is a soft semi-neighborhood of each point of $\bigcap_{\alpha \in E} F(\alpha)$.

**Proof.** The proof directly follows form the definition of soft semi-neighborhood.

**Proposition 3.4.** Let $(X, \tau, E)$ be a soft topological space over X. Then

(1) each $x \in X$ has a soft semi-neighborhood;

(2) if $(F,E)$ and $(G,E)$ are soft semi-neighborhoods of some $x \in X$, then $(F,E) \cap (G,E)$ is also a soft semi-neighborhood of $x$.

(3) if $(F,E)$ is a soft semi-neighborhood of $x \in X$ and $(F,E) \widetilde{\subset} (G,E)$, then $(G,E)$ is also a soft semi-neighborhood of $x \in X$.

**Proof.** The proof is easy and thus omitted.

In [3], we have the following properties of soft semi-interior:

**Theorem 3.5.** Let $(X, \tau, E)$ be a soft topological space over $X$ and $(F,E)$ and $(G,E)$ are soft sets over $X$. Then

(1) $int^s \Phi = \Phi$ and $int^s \widetilde{X} = \widetilde{X}$.

(2) $int^s(F,E) \widetilde{\subset} (F,E)$.

(3) $(int^s(int^s(F,E))) = (F,E)$

(4) $(F,E)$ is a soft semi-open set if and only if $int^s(F,E) = (F,E)$

(5) $(F,E) \widetilde{\subset} (G,E)$ implies $int^s(F,E) \widetilde{\subset} int^s(G,E)$

(6) $int^s(F,E) \cup int^s(G,E) \widetilde{\subset} int^s((F,E) \cup (G,E))$.

Now we prove the following:

**Theorem 3.6.** Let $(X, \tau, E)$ be a soft topological space over $X$ and $(F,E)$ and $(G,E)$ are soft sets over $X$. Then

$$int^s(F,E) \cap int^s(G,E) = int^s((F,E) \cap (G,E))$$

**Proof.** Using (5) of above theorem, we have $((F,E) \cap (G,E)) \widetilde{\subset} (F,E)$, $((F,E) \cap (G,E)) \widetilde{\subset} (G,E)$ implies $int^s((F,E) \cap (G,E)) \widetilde{\subset} int^s(F,E)$, $int^s((F,E) \cap (G,E)) \widetilde{\subset} int^s(G,E)$ so that $int^s((F,E) \cap (G,E)) \widetilde{\subset} int^s(F,E) \cap int^s(G,E)$. Also, since $int^s(F,E) \widetilde{\subset} (F,E)$ and $int^s(G,E) \widetilde{\subset} (G,E)$ implies $int^s(F,E) \cap int^s(G,E) \widetilde{\subset} ((F,E) \cap (G,E))$ so that $int^s(F,E) \cap int^s(G,E)$ is a soft semi-open subset of $((F,E) \cap (G,E))$. Hence $int^s(F,E) \cap int^s(G,E) \widetilde{\subset} int^s((F,E) \cap (G,E))$. Thus $int^s(F,E) \cap int^s(G,E) = int^s((F,E) \cap (G,E))$.

**Definition 3.7.** Let $(X, \tau, E)$ be a soft topological space over $X$ then soft semi-exterior of soft set $(F,E)$ over $X$ is denoted by $ext^s(F,E)$ and is defined as $ext^s(F,E) = int^s((F,E)')$.

Thus $x$ is called a soft semi-exterior point of $(F,E)$ if there exists a soft semi-open set $(G,E)$ such that $x \in (G,E) \widetilde{\subset} (F,E)'$. We observe that $ext^s(F,E)$ is the largest soft semi-open set contained in $(F,E)'$.

**Example 3.8.** Let $X = \{h_1, h_2, h_3\}$, $E = \{e_1, e_2\}$ and

$\tau = \{\Phi, \widetilde{X}, (G_1,E), (G_2,E), (G_3,E), (G_4,E), (G_5,E), (G_6,E), (G_7,E)\}$

where $(G_1,E)$, $(G_2,E)$, $(G_3,E)$, $(G_4,E)$, $(G_5,E)$, $(G_6,E)$ and $(G_7,E)$ are soft sets over $X$, defined as follows:

$$G_1(e_1) = \{h_1, h_2\}, \ G_1(e_2) = \{h_1, h_2\},$$
$$G_2(e_1) = \{h_2\}, \ G_2(e_2) = \{h_1, h_3\},$$
$$G_3(e_1) = \{h_2, h_3\}, \ G_3(e_2) = \{h_1\},$$
$$G_4(e_1) = \{h_2\}, \ G_4(e_2) = \{h_1\},$$
$$G_5(e_1) = \{h_1, h_2\}, \ G_5(e_2) = X,$$
$$G_6(e_1) = X, \ G_6(e_2) = \{h_1, h_2\},$$
$$G_7(e_1) = \{h_2, h_3\}, \ G_7(e_2) = \{h_1, h_3\}.$$

Then $\tau$ defines a soft topology on $X$ and hence $(X, \tau, E)$ is a soft topological space over $X$. Let us take $(F, E) = \{\{h_1\}, \{h_3\}\}$ be soft set of soft topological space $X$. Then $(F, E)' = \{\{h_2, h_3\}, \{h_1, h_2\}\}$ and so, $ext^s(F, E) = int^s((F, E)') = \{\{h_2, h_3\}, \{h_1, h_2\}\}$, since $(G, E) = \{\{h_2, h_3\}, \{h_1, h_2\}\}$ is soft semi-open set.

The proof of the following theorem directly follows from definition of soft semi-exterior:

**Theorem 3.9.** Let $(F, E)$ be a soft set of soft topological space over X. Then

$$(1) ext^s(F, E) = int^s((F, E)').$$
$$(2) ext^s((F, E) \cup (G, E)) = ext^s(F, E) \cap ext^s(G, E).$$
$$(3) ext^s(F, E) \cup ext^s(G, E) \widetilde{\subset} ext^s((F, E) \cap (G, E)).$$

The following example shows that the equality does not holds in Theorem 3.9(3).

**Example 3.10.** In Example 3.8, let us take e $(F, E) = \{\{h_3\}, \{h_3\}\}$, $(G, E) = \{\{h_1, h_2\}, \{h_1, h_2\}\}$. Now

$$ext^s((F, E) \cap (G, E)) = int^s(((F, E) \cap (G, E))') = int^s((\{\{h_3\}, \{h_3\}\} \cap \{\{h_1, h_2\}, \{\{h_1, h_2\}\}\})')$$
$$= int^s((\tilde{\phi})') = int^s(\tilde{X}) = \tilde{X}.$$

Also $ext^s(F, E) = int^s((F, E)') = int^s\{\{h_1, h_2\}, \{h_1, h_2\}\} = \{\{h_1, h_2\}, \{h_1, h_2\}\}$

and $ext^s(G, E) = int^s((G, E)') = int^s(\{\{h_1, h_2\}, \{h_1, h_2\}\}') = int^s\{\{h_3\}, \{h_3\}\} = \tilde{\phi}$. Thus $ext^s(F, E) \cup ext^s(G, E) = \{\{h_1, h_2\}, \{h_1, h_2\}\}$.

**Definition 3.11.** Let $(X, \tau, E)$ be a soft topological space over $X$ then soft semi-boundary of soft set $(F, E)$ over $X$ is denoted by $bd^s(F, E)$ and is defined as

$$bd^s(F, E) = (int^s(F, E) \cup ext^s(F, E))'.$$

**Remark 3.12.** From the above definition it follows directly that the soft sets $(F, E)$ and $(F, E)'$ have same soft semi-boundary.

**Theorem 3.13.** Let $(F, E)$ be a soft set of soft topological space over X. Then the following hold:

$(1) cl^s(F, E) = int^s(F, E) \cup bd^s(F, E).$

$(2) bd^s(F, E) = cl^s(F, E) \cap cl^s(F, E)' = cl^s(F, E) - int^s(F, E).$

$(3) (bd^s(F, E))' = int^s(F, E) \cup int^s((F, E)') = int^s(F, E) \cup ext^s(F, E).$

$(4) int^s(F, E) = (F, E) \setminus bd^s(F, E).$

**Proof.**

(1). $int^s(F,E) \cup bd^s(F,E) = int^s(F,E) \cup (cl^s(F,E) \cap cl^s(F,E)')$
$= [int^s(F,E) \cup cl^s(F,E)] \cap [int^s(F,E) \cup cl^s(F,E)']$
$= cl^s(F,E) \cap [int^s(F,E) \cup (int^s(F,E))']$
$= cl^s(F,E) \cap (int^s(F,E) \cup (int^s(F,E))')$
$= cl^s(F,E) \cap \widetilde{X}$
$= cl^s(F,E).$

(2). $bd^s(F,E) = cl^s(F,E) - int^s(F,E)$
$= cl^s(F,E) \cap (int^s(F,E))'$
$= cl^s(F,E) \cap cl^s(F,E)'$ (by Theorem 3.8(2)[3]).

(3). $int^s(F,E) \cup int^s((F,E)') = ((int^s(F,E))')' \cup (int^s(((F,E)')))'$
$= [(int^s(F,E))' \cap int^s((F,E)')]'$
$= [cl^s(F,E)' \cap cl^s(F,E)]'$
$= (bd^s(F,E))'.$

(4). $(F,E) \setminus bd^s(F,E) = (F,E) \cap bd^s(F,E)'$
$= (F,E) \cap (int^s(F,E) \cup int^s((F,E)'))$ (by (1))
$= [(F,E) \cap int^s(F,E)] \cup [(F,E) \cap int^s((F,E)')]$
$= int^s(F,E) \cup \Phi$
$= int^s(F,E).$

**Remark 3.14.** By (2) of above Theorem 3.13, it is clear that
(i) $bd^s(F,E)$ is a smallest soft semi-closed set over $X$ containing $(F,E)$
(ii) $bd^s(F,E) = bd^s(F,E)'$.

**Theorem 3.15.** Let $(F,E)$ be soft set of soft topological space over X. Then:
(1) $(F,E)$ is soft semi-open set over $X$ if and only if $(F,E) \cap bd^s(F,E) = \Phi$.
(2) $(F,E)$ is soft semi-closed set over $X$ if and only if $bd^s(F,E) \widetilde{\subset} (F,E)$.

**Proof.** (1). Let $(F,E)$ be a soft semi-open set over $X$. Then $int^s(F,E) = (F,E)$ implies

$$(F,E) \cap bd^s(F,E) = int^s(F,E) \cap bd^s(F,E)$$
$$= \Phi.$$

Conversely, let $(F,E) \cap bd^s(F,E) = \Phi$. Then $(F,E) \cap cl^s(F,E) \cap cl^s((F,E)') = \Phi$ or $(F,E) \cap cl^s((F,E)') = \Phi$ or $cl^s((F,E)') \tilde{\subset} (F,E)'$, which implies $(F,E)'$ is soft semi-closed and hence $(F,E)$ is soft semi-open set.

(2). Let $(E,E)$ be soft semi-closed set over $X$. Then $cl^s(F,E) = (F,E)$. Now $bd^s(F,E) = cl^s(F,E) \cap cl^s((F,E)') \tilde{\subset} cl^s(F,E) = (F,E)$. That is, $bd^s(F,E) \tilde{\subset} (F,E)$.

Conversely, $bd^s(F,E) \tilde{\subset} (F,E)$. Then $bd^s(F,E) \cap (F,E)' = \Phi$. Since $bd^s(F,E) = bd^s((F,E)') = \Phi$, we have $bd^s((F,E)') \cap (F,E)' = \Phi$. By (1), $(F,E)'$ is soft semi-open and hence $(F,E)$ is soft semi-closed.

**Theorem 3.16.** Let $(F,E)$ be a soft set of soft topological space over X. Then the following hold:

$$(1) bd^s(F,E) \cap int^s(F,E) = \Phi.$$
$$(2) cl^s int^s(F,E) = (F,E) \setminus bd^s(F,E).$$

**Proof.** (1) follows form Theorem 3.13(2) and (2) follows directly by the definition of the soft boundary.

**Theorem 3.17.** Let $(F,A)$ and $(G,B)$ be soft sets of soft topological space over X. Then the following hold:

$(1) bd^s((F,A) \cup (G,B)) \tilde{\subset} [bd^s(F,A) \cap ((G,B)')] \cup [bd^s(G,B) \cap cl^s((F,A)')]$.

$(2) bd^s[(F,A) \cap (G,B)] \tilde{\subset} [bd^s(F,A) \cap cl^s(G,B)] \cup [bd^s(G,B) \cap cl^s(F,A)]$.

**Proof.**

(1).
$$\begin{aligned}
bd^s((F,A) \cup (G,B)) &= cl^s((F,A) \cup (G,B)) \cap cl^s(((F,A) \cup (G,B))') \\
&= (cl^s(F,A) \cup cl^s(G,B)) \cap cl^s((F,A)' \cap (G,B)') \\
&\tilde{\subset} (cl^s(F,A) \cup cl^s(G,B)) \cap [cl^s(F,A)' \cap cl^s(G,B)'] \\
&= (cl^s(F,A) \cap cl^s(F,A)') \cap (cl^s(G,B)' \cup cl^s(G,B)) \cap [cl^s(F,A)' \cap cl^s(G,B)'] \\
&= (bd^s(F,A) \cap cl^s(G,B)') \cup (bd^s(G,B) \cap cl^s(F,A)') \\
&\tilde{\subset} bd^s(F,A) \cup bd^s(G,B).
\end{aligned}$$

(2).
$$\begin{aligned}
bd^s[(F,A) \cap (G,B)] &= cl^s((F,A) \cap (G,B)) \cap cl^s((F,A) \cap (G,B))' \\
&\tilde{\subset} [cl^s(F,A) \cap cl^s(G,B)] \cap [cl^s((F,A)' \cup (G,B)')] \\
&= [cl^s(F,A) \cap cl^s(G,B)] \cap [cl^s(F,A)' \cup cl^s(G,B)'] \\
&= [(cl^s(F,A) \cap cl^s(G,B)) \cap cl^s(F,A)'] \cup [(cl^s(F,A) \cap cl^s(G,B)) \cap cl^s(G,B)'] \\
&= (bd^s(F,A) \cap cl^s(G,B)) \cup (cl^s(F,A) \cap bd^s(G,B)).
\end{aligned}$$

**Theorem 3.18.** Let $(F,E)$ be a soft set of soft topological space over X. Then the following holds:
$$bd^s(bd^s(bd^s(F,E))) = bd^s(bd^s(F,E)).$$

**Proof.**

$$bd^s(bd^s(bd^s(F,E))) = cl^s(bd^s(bd^s(F,E))) \cap cl^s((bd^s(bd^s(F,E)))')$$
$$= (bd^s(bd^s(F,E))) \cap cl^s((bd^s(bd^s(F,E)))') \ldots \quad (1)$$

Now consider

$$((bd^s(bd^s(F,E)))') = [cl^s(bd^s(F,E)) \cap ((bd^s(F,E))')]'$$
$$= [bd^s(F,E) \cap cl^s(bd^s(F,E))']'$$
$$= (bd^s(F,E))' \cup (cl^s(bd^s(F,E))')'$$

Therefore

$$cl^s((bd^s(bd^s(F,E)))') = cl^s[cl^s((bd^s(F,E))') \cup (cl^s((bd^s(F,E))'))']$$
$$= cl^s(cl^s((bd^s(F,E))')) \cup cl^s((cl^s((bd^s(F,E))'))')$$
$$= (G,E) \cup ((cl^s((bd^s(G,E))'))') = \tilde{X} \ldots \quad (2)$$

Where $(G,E) = cl^s(cl^s((bd^s(F,E))'))$. From (1) and (2), we have

$$bd^s(bd^s(bd^s(F,E))) = bd^s(bd^s(F,E)) \cap \tilde{X} = bd^s(bd^s(F,E)).$$

**Theorem 3.19.** Let $(F,E)$ and $(G,E)$ be a soft open sets of soft topological space over X. Then the following hold:

(1) $((F,E) \setminus int^s(G,E)) \tilde{\subset} (F,E)^\circ \setminus int^s(G,E)$

(2) $bd^s int^s(F,E) \tilde{\subset} bd^s(F,E)$

**Proof.**

(1) $((F,E) \setminus int^s(G,E)) = ((F,E) \cap int^s(G,E)')$
$$= int^s(F,E) \cap int^s((G,E)')$$
$$= int^s(F,E) \cap (cl^s(G,E))' \qquad (by\ Theorem\ 3.8(1)[3])$$
$$= int^s(F,E) \setminus cl^s(G,E)$$
$$\tilde{\subset} (F,E)^\circ \setminus int^s(G,E).$$

(2) $bd^s int^s(F,E) = cl^s int^s(F,E) \cap cl^s((int^s(F,E))')$
$$\tilde{\subset} cl^s int^s(F,E) \cap cl^s(cl^s((F,E)')) \qquad (by\ Theorem\ 3.8(2)[3])$$
$$\tilde{\subset} cl^s(F,E) \cap cl^s((F,E)') = bd^s(F,E).$$

**Theorem 3.20.** Let $(F,E)$ be a soft set of soft topological space over X. Then $bd^s(F,E) = \Phi$ if and only if $(F,E)$ is a soft semi-closed set and a soft semi-open set.

**Proof.** Suppose that $bd^s(F,E) = \Phi$.

**(i)** First we prove that $(F,E)$ is a soft semi-closed set. Consider

$$bd^s(F,E) = \Phi \Rightarrow cl^s(F,E) \cap cl^s((F,E)') = \Phi$$
$$\Rightarrow cl^s(F,E) \tilde{\subset} (cl^s((F,E)'))' = int^s(F,E) \qquad (by\ Theorem\ 3.8(3)[3])$$
$$\Rightarrow cl^s(F,E) \tilde{\subset} (F,E) \Rightarrow cl^s(F,E) = (F,E)$$

This implies that $(F,E)$ is a soft semi-closed set.

**(ii)** Using (i), we now prove that $(F,E)$ is a soft semi-open set. $bd^s(F,E) = \Phi \Rightarrow$

$cl^s(F,E) \cap cl^s((F,E)')$ or $(F,E) \cap int^s((F,E))' = \Phi \Rightarrow (F,E) \tilde{\subseteq} int^s(F,E) \Rightarrow int^s(F,E) = (F,E)$. This implies that $(F,E)$ is a soft semi-open set.

Conversely, suppose that $(F,E)$ is soft semi-open and soft semi-closed set. Then
$$\begin{aligned} bd^s(F,E) &= cl^s(F,E) \cap cl^s((F,E)') \\ &= cl^s(F,E) \cap (int^s(F,E))' \quad (by\ Theorem\ 3(4) \\ &= (F,E) \cap (F,E)' = \Phi. \end{aligned}$$
This completes the proof.

**Definition 3.21.** Let $(X,\tau,E)$ be a soft topological space over $X$ and $x \in (F,E)$. If there is a soft semi-open set $(G,E)$ such that $x \in (G,E)$, then $(G,E)$ is called a soft semi-open neighborhood (or soft semi-open nbd) of $x$. The set of all soft semi-open nbds of $x$, denoted $\tilde{n}bd^s(x)$, is called the soft semi-open nbd systems of $x$; that is,
$$nbd^s(x) = \{(G,E) : (G,E) \in S.S.O(X) : x \in (G,E)\}.$$

The following theorem gives important properties of soft semi-open nbd system:

**Proposition 3.22.** Let $(X,\tau,E)$ be a soft topological space over $X$ and $(G,E),(H,E) \tilde{\subseteq} (F,E)$. Then the collection of soft semi-open nbd $\tilde{n}bd^s(x)$ at $x$ in $(X,\tau,E)$ has the following properties:

(1) If $(G,E) \in \tilde{n}bd^s(x)$, then $x \in (G,E)$.

(2) If $(G,E),(H,E) \in \tilde{n}bd^s(x)$, then $(G,E) \tilde{\cap} (H,E) \in \tilde{n}bd^s(x)$.

(3) If $(G,E) \in \tilde{n}bd^s(x)$ and $(G,E) \tilde{\subseteq} (H,E)$, then $(H,E) \in \tilde{n}bd^s(x)$.

(4) If $(G,E) \in \tilde{n}bd^s(x)$, then there is an $(H,E) \in \tilde{n}bd^s(x)$ such that $(G,E) \in \tilde{n}bd^s(y)$, for each $y \in (H,E)$.

(5) $(G,E) \tilde{\subseteq} (F,E)$ is soft open if and only if $(G,E)$ contains a soft semi-open nbd of each of its points.

**Proof.** (1) is obvious, since $(G,E)$ is a soft semi-open nbd of $x \in (F,E)$. So $(G,E)$ is a soft semi-open set such that $x \in (G,E)$.

(2). If $(G,E),(H,E) \in \tilde{n}bd^s(x)$, then there exist soft semi-open sets $(I,E)$ and $(J,E)$ such that $x \in (I,E) \tilde{\subseteq} (G,E)$ and $x \in (J,E) \tilde{\subseteq} (H,E)$. Therefore $x \in (I,E) \tilde{\cap} (J,E) \tilde{\subseteq} (G,E) \tilde{\cap} (H,E)$ and hence $(G,E) \tilde{\cap} (H,E) \in \tilde{n}bd^s(x)$.

(3). Since $(G,E) \in \tilde{n}bd^s(x)$, therefore there exists a soft semi-open set $(I,E)$ such that $x \in (I,E) \tilde{\subseteq} (G,E)$. Then $x \in (I,E) \tilde{\subseteq} (G,E) \tilde{\subseteq} (H,E)$ or $x \in (I,E) \tilde{\subseteq} (H,E)$. Hence $(H,E) \in \tilde{n}bd^s(x)$.

(4) Since $(G,E) \in \tilde{n}bd^s(x)$, then $x \in (H,E) \tilde{\subseteq} (G,E)$, for $(H,E)$ soft semi-open in $(F,E)$. Since $x \in (H,E) \tilde{\subseteq} (H,E)$ then $(H,E) \in \tilde{n}bd^s(x)$. If $y \in (H,E)$, then by (3) $(H,E) \tilde{\subseteq} (G,E)$ implies $(G,E) \in \tilde{n}bd^s(x)$, for each $y \in (H,E)$.

(5). (i) Suppose $(G,E)$ is a soft semi-open in $(F,E)$, then $x \in (G,E) \tilde{\subseteq} (G,E)$ implies $(G,E)$ is a soft semi-open nbd of each $x \in (G,E)$.

(ii) If each $x \in (G,E)$ has a soft semi-open nbd $(H,E)_x \tilde{\subseteq} (G,E)$, then $(G,E) = \{x : x \in (G,E)\} \tilde{\subseteq} \bigcup_{x \in (G,E)} (H,E)_x \tilde{\subseteq} (G,E)$ or $(G,E) = \bigcup_{x \in (G,E)} (H,E)_x$. This gives $(G,E)$ is soft semi-open in $(F,E)$.

**Definition 3.23.** Let $(X, \tau, E)$ be a soft topological space over $X$. A soft semi-nbd base at $x \in (F, E)$ is a subcollection $\tilde{s}nbd^s(x)$ of soft semi-nbd $\tilde{n}bd^s(x)$ having the property that each $(G, E) \in \tilde{n}bd^s(x)$ contains some $(H, E) \in \tilde{s}nbd^s(x)$. That is, $\tilde{n}bd^s(x)$ must be determined by $\tilde{s}nbd^s(x)$ as follows:

$$\tilde{n}bd^s(x) = \{(G, E) \tilde{\subseteq} (F, E) : (H, E) \tilde{\subseteq} (G, E), \text{ for some } (H, E) \in \tilde{s}nbd^s(x)\}.$$

Each $(H, E) \in \tilde{s}nbd^s(x)$ is called a basic soft semi-open neighborhood of $x$.

For the soft basic semi-nbd system, the following properties directly follows by the corresponding properties of Proposition 3.22.:

**Proposition 3.24.** Let $(X, \tau, E)$ be a soft topological space over $X$ and for each $x \in (F, E)$, let $\tilde{s}nbd^s(x)$ be a soft semi-nbd base at $x$. Then

(1) If $(H, E) \in \tilde{s}nbd^s(x)$, then $x \in (H, E)$.

(2) If $(F_1, E), (F_2, E) \in \tilde{s}nbd^s(x)$, then there is some $(F_3, E) \in \tilde{s}nbd^s(x)$ such that $(F_3, E) \tilde{\subseteq} (F_1, E) \tilde{\cap} (F_2, E)$.

(3) If $(H, E) \in \tilde{s}nbd^s(x)$, then there is some $(F_0, E) \in \tilde{s}nbd^s(x)$ such that if $z \in (F_0, E)$, then there is some $(I, E) \in \tilde{s}nbd^s(z)$ with $(I, E) \tilde{\subseteq} (H, E)$.

(4) $(G, E) \tilde{\subseteq} (F, E)$ is soft semi-open if and only if $(G, E)$ contains a soft semi-basic nbd of each of its points.

**Conclusion:** Soft set theory is very important during the study towards possible applications in classical and non classical logic. In recent years many researchers worked on the findings of structures of soft sets theory initiated by Molodtsov and applied to many problems having uncertainties. It is worth mentioning that soft topological spaces based on soft set theory which is a collection of information granules is the mathematical formulation of approximate reasoning about information systems. In the present work, we continued to investigate the properties of soft semi-open sets and soft semi-closed sets in soft topological spaces. We defined soft semi-exterior, soft semi-boundary, soft semi-open neighborhood and soft semi-open neighborhood systems in soft topological spaces. Moreover we discussed the characterizations of soft semi-closed and soft semi-open sets via soft semi-interior, soft semi-exterior, soft semi-closure and soft semi-boundary and have established several interesting properties. We hope that the findings in this paper will help the researchers to enhance and promote the further study on soft topology to carry out general framework for the applications in practical life.